\documentclass[10pt,twocolumn]{article}

\usepackage{graphicx}
\usepackage{color}
\usepackage[utf8]{inputenc}
\usepackage{enumerate}
\usepackage{amsmath}
\usepackage{amsthm}
\usepackage{amssymb}
\usepackage{accents}
\usepackage{subfig}
\usepackage[style=ieee,backend=bibtex]{biblatex}
\usepackage{algorithm}
\usepackage{algorithmic}
\usepackage{hyperref}
\usepackage{nameref}
\usepackage{eurosym}
\usepackage{flushend}
\usepackage{setspace}
\hypersetup{
     colorlinks=true, breaklinks=true, bookmarksopen=true,
     urlcolor=black,  citecolor=black, linkcolor=black
}      

\usepackage{framed}

\usepackage{todonotes}
\usepackage{bm}
\newcommand{\bl}{ \bm{\lambda}} 
\newcommand{\bg}{ \bm{\gamma}} 
\newcommand{\ba}{ \bm{\alpha}} 
\newcommand{\bb}{ \bm{\beta}} 
\newcommand{\tinyt}{{\scriptscriptstyle(t)}}

\usepackage{tikz}
\usepackage[europeanresistors]{circuitikz}

\newcommand{\image}[4][8.89cm] 
{
\begin{figure}[!ht]
\begin{center}
	\includegraphics[width=#1]{#2}
	\caption{#3}
	\label{fig:#4}
\end{center}
\end{figure}
}

\makeatletter
\def\section{\@startsection{section}{1}{\z@}{1.5ex plus 1.5ex minus 0.5ex}{0.7ex plus 1ex minus 0ex}{\normalfont\normalsize\centering\scshape}}%
\def\subsection{\@startsection{subsection}{2}{\z@}{1.5ex plus 1.5ex minus 0.5ex}{0.7ex plus .5ex minus 0ex}{\normalfont\normalsize\itshape}}%
\makeatother

\title{\vspace{-2em}{\Huge Relaxations for multi-period optimal power flow \\[0.3em] problems with discrete decision variables}}
\author{Q. Gemine, D. Ernst, Q. Louveaux, B. Corn\'elusse\\
 Department of Electrical Engineering  and Computer Science,\\
        University of Li\`ege,  Li\`ege 4000, Belgium\\
\{qgemine,dernst,q.louveaux,bertrand.cornelusse\}@ulg.ac.be}
\date{}

\allowdisplaybreaks

\begin{document}
\maketitle
\setstretch{0.8}
\begin{Abstract}
\textbf{\small We consider a class of optimal power flow (OPF) applications where some loads offer a modulation service in exchange for an activation fee. These applications can be modeled as multi-period formulations of the OPF with discrete variables that define mixed-integer non-convex mathematical programs. We propose two types of relaxations to tackle these problems. One is based on a Lagrangian relaxation and the other is based on a network flow relaxation. Both relaxations are tested on several benchmarks and, although they provide a comparable dual bound, it appears that the constraints in the solutions derived from the network flow relaxation are significantly less violated.}
\end{Abstract}


\setstretch{1}

\section{Introduction}
\label{sec:introduction}

Many power system applications that require solving an optimal power flow (OPF) problem share two features. 
 Firstly, these applications are multi-period because of the evolution of market prices, of the ramping limits of generation units and of the behavior of static and flexible loads. 
 Secondly they contain integer decision variables to model the acceptance or the rejection of bids, or the start up of some generation units.
 As a first example, the day-ahead energy market in Europe computes spot prices based on supply and demand offers. This application has a multi-period and discrete nature because of the ``block bids'', and because of some ramping constraints. Active power flows are constrained by a simple network flow model. Operational constraints on reactive power, voltage and current are aggregated in the arc capacities of the network flows. More realistic (so called ``flow based''~\cite{aguado2012flow}) network models are emerging, but they are still a linear approximation of the set of feasible flows around a foreseen operation point. 
As a second example, new applications arising in distribution networks such as operational planning aim at avoiding  the congestion of network elements and minimizing the curtailment of renewable energy sources. To benefit from the flexibility of customers, it is necessary to account for the time-coupled nature of the problem, and integer variables can be used to model the reservation of that flexibility.  The physical characteristics of the network are different from those of transmission systems and DC power flow approximations can hardly be used. 

Hence depending on the complexity of the primary goal of the application and its scale, it is often mandatory to resort to a relaxation of the non-convex network constraints so as to devise a robust and fast algorithm. Also, a common characteristic of these applications is that the main decision variables are the power injections, and especially active power flows as they underlie most of the financial transactions. The other variables (voltage, current) can be viewed as a consequence of the power flows in the network, and we must ensure that these consequences stay within the operational limits. 
These observations motivate the relaxation algorithms studied in this paper. We focus on relaxations that decompose the problem into one subproblem that works exclusively with active and reactive power flows but encompasses the multi-period and discrete aspects, and subproblems that assert that for each time step those flows do not violate voltage and other technical limits. 
After the precise statement of the discrete multi-period optimal power flow we are targeting in Section~\ref{sec:problem-statement} and a review of the recent literature on these topics in Section~\ref{sec:liter-revi-single}, we propose two relaxations achieving these goals in Section~\ref{sec:relax}. 
The first relaxation is a straightforward generalization of the Lagrangian relaxation (LR) of \cite{phan2012} to this problem. 
The downside of this LR scheme is that the power related subproblem lacks information on the network topology.
The second relaxation builds on a network flow reformulation of the original problem by introducing link-flow variables. It is then relaxed into a convex problem by substituting non-linear terms with their convex envelopes. Small semidefinite programming (SDP) relaxations are used to translate operational limits into bounds of voltage and link-flow variables.
Section~\ref{sec:results} compares the two proposed approaches on several test systems, whereas Section~\ref{sec:conclusion} concludes and gives directions of further research.

\section{General problem statement}
\label{sec:problem-statement}

We consider the problem of finding the optimal operation of a set $\mathcal{D}$ of devices (i.e. loads and generators) over a certain time horizon while maintaining the network and the devices within operational limits. The network is defined as a set $\mathcal{L}$ of links, that is lines, cables or transformers that define pairwise connections between elements of the set $\mathcal{B}$ of buses. Several devices can be connected to a single bus. The time horizon is modeled by a set $\mathcal{T}$ of periods.
We denote by $\mathcal{F} \subset \mathcal{D}$ the flexible loads. The consumption of a flexible load can be modulated around a baseline profile. In particular, we use the flexibility model presented in \cite{gemineactive}, where the right to modulate a flexible load is conditioned to the payment of an availability fee. The operational constraints associated to these loads are upward and downward modulation limits as well as an energy constraint, stating that any modulation should consume the same amount of energy than the baseline profile.
We use the following notations throughout this paper, where the superscript $(t)$ refers to period $t$:
\begin{itemize}
\item $\bm{P}^{(t)} \in \mathbb{R}^{|\mathcal{D}|}$, the active power injections of devices (positive when power flows from the device to the network);
\item $\bm{Q}^{(t)} \in \mathbb{R}^{|\mathcal{D}|}$, the reactive power injections of devices (same sign convention as $\bm{P}^{(t)}$);
\item $\bm{d} \in \{0,1\}^{|\mathcal{F}|}$, the availability indicators of flexible loads;
\item $\bm{c}_f \in \mathbb{R}_+^{|\mathcal{F}|}$, the availability costs of the flexible loads;
\item ${\bm{P}^{bl}}^{(t)} \in \mathbb{R}^{|\mathcal{F}|}$, the active power injections of flexible loads when operating at their baseline;
\item $\overline{\bm{P}}^{(t)}$ and $\underline{\bm{P}}^{(t)} \in \mathbb{R}^{|\mathcal{B}|}$, the bounds on active power injection of devices;
\item $\mathbb{A} \in \mathbb{R}^{N_c \times 2|\mathcal{D}|}$ and $\bm{a} \in \mathbb{R}^{N_c}$, matrix and vector modeling the P-Q capability of the devices (with $N_c$ the total number of linear constraints between $\bm{P}^{(t)}$ and $\bm{Q}^{(t)}$);
\item $\mathbb{M} \in \{0,1\}^{\mathcal{B} \times \mathcal{D}}$, mapping from devices to
  buses ($M_{i,j} = 1$ if device $j$ is connected to bus \nolinebreak
  $i$ and $0$ otherwise);
\item $\bm{e}^{(t)} \in \mathbb{R}^{|\mathcal{B}|}$, the real part of the voltage at buses;
\item $\bm{f}^{(t)} \in \mathbb{R}^{|\mathcal{B}|}$, the imaginary part of the voltage at buses;
\item $\overline{\bm{V}}$ and $\underline{\bm{V}} \in \mathbb{R}^{|\mathcal{B}|}$, the limits on the voltage magnitudes;
\item $g_{ij}$ the conductance of link $(i,j) \in \mathcal{L}$;
\item $b_{ij}$ the susceptance of link $(i,j) \in \mathcal{L}$.
\end{itemize}
The decision variables are the subset of the active and reactive power injections for which the bounds $\overline{P}^{(t)}_k$ and $\underline{P}^{(t)}_k$ are not equal ($k \in \mathcal{D})$, the voltage at all buses, and the discrete decision variables $\bm{d}$. The notion of optimal operation is defined by a generic cost function $f(\bm{P})$ (linear or a convex quadratic) that we want to minimize together with ${\bm{c}_f\cdot\bm{d}}$, the availability fees of flexible loads.
The whole problem is modeled in (\ref{eq:1}-\ref{eq:10}) where we use the notation $\bm{P}$, $\bm{Q}$, $\bm{e}$ and $\bm{f}$ to denote the concatenation of, respectively, the vectors $\bm{P}_t$, $\bm{Q}_t$, $\bm{e}_t$ and $\bm{f}_t$ for all $t \in \mathcal{T}$.
\begin{align}
\min_{\substack{\bm{P}, \bm{Q} \\ \bm{d}, \bm{e}, \bm{f}}} & & & f(\bm{P}) + \bm{c}_f \cdot \bm{d}\label{eq:1}\\
\text{s.t.}& & & \mathbf{d}  \in \{0,1\}^{|\mathcal{F}|}\label{eq:2}\\
& & & \hspace{-0.75em}\forall~t \in \mathcal{T}:\nonumber\\
& & & \underline{\bm{P}}^{(t)} \leq \bm{P}^{(t)} \leq \overline{\bm{P}}^{(t)}\label{eq:3}\\
& & & \mathbb{A}\begin{pmatrix}\bm{P}^{(t)}\\\bm{Q}^{(t)}\end{pmatrix} \leq \bm{a}\label{eq:4}\\
& & & \hspace{-0.75em}\forall k \in \mathcal{F}:\nonumber\\
& & & \sum_{t \in \mathcal{T}} \Big(P_k^{(t)} - {P_k^{bl}}^{(t)}\Big) = 0 \label{eq:5}\\
& & & \hspace{-0.75em}\forall (t,k) \in \mathcal{T} \times \mathcal{F}:\nonumber\\
& & & P_k^{(t)} \geq (1 - d_k){P_k^{bl}}^{(t)} + d_k \underline{P}_k^{(t)}\label{eq:6}\\
& & & P_k^{(t)} \leq (1-d_k){P_k^{bl}}^{(t)} + d_k \overline{P}_k^{(t)}\label{eq:7}\\
& & & \hspace{-0.75em}\forall (t,i) \in \mathcal{T} \times \mathcal{B}:\nonumber\\
& & &(\mathbb{M}~\bm{P}^{(t)})_i = \sum_{j \in \mathcal{N}(i)} \Big(g_{ij} ({e_i^\tinyt}^2 + {f_i^\tinyt}^2 - e_i^\tinyt e_j^\tinyt \nonumber\\ 
& & & \hspace{4.4em} - f_i^\tinyt f_j^\tinyt) + b_{ij} (e_i^\tinyt f_j^\tinyt - f_i^\tinyt e_j^\tinyt)\Big)\label{eq:8}\\
& & &(\mathbb{M}~\bm{Q}^{(t)})_i = \sum_{j \in \mathcal{N}(i)} \Big(b_{ij} (e_i^\tinyt e_j^\tinyt + f_i^\tinyt f_j^\tinyt- {e_i^\tinyt}^2\nonumber \\ 
& & & \hspace{4.4em} - {f_i^\tinyt}^2) + g_{ij} (e_i^\tinyt f_j^\tinyt - f_i^\tinyt e_j^\tinyt)\Big)\label{eq:9}\\
& & &\underline{V}_i^2 \leq {e_i^\tinyt}^2 + {f_i^\tinyt}^2 \leq \overline{V}_i^2\label{eq:10}
\end{align}
This is a mixed-integer non-convex mathematical program where the non-convexity comes from constraints (\ref{eq:8}-\ref{eq:10}). In addition, the electrical variables (i.e. powers and voltages) are coupled over the set $\mathcal{T}$ of periods because of the time-coupling constraints (\ref{eq:5}-\ref{eq:7}) that model the flexible loads.

\section{Literature review}
\label{sec:liter-revi-single}

We first review the methods designed to solve "static" OPF problems, in the sense that the problem has no temporal aspect. We then review the literature on multi-period OPF, which is a scale up of a static problem caused by time coupling constraints on power variables. Finally, we review the literature on works where some discrete variables have been introduced in the OPF problem to model the ability to act on power injections or withdrawals, that is, problems comparable to the problem introduced in Section~\ref{sec:problem-statement}.

Optimal power flow problems, although non-convex, have been for long solved using local non-linear optimization methods. Interior-point methods are probably the most widespread class of methods dedicated to this problem \cite{capitanescu2007ipm}. If the solution they provide has no guarantee to be globally optimal, they have been made popular by their convergence speed and their ability to solve fairly efficiently problems of large dimension.

Recently, SDP was successfully applied as a convex relaxation to the OPF problem \cite{lavaei2012}. The OPF is formulated over all the degree 2 monomials of the real and imaginary parts of the voltage variables. Dropping the rank 1 constraint of the corresponding matrix yields the SDP relaxation. For technical reasons, the dual of this SDP relaxation is solved (strong duality holds). When the duality gap is zero, a primal feasible optimal solution to the original OPF problem can be recovered from the solution of the dual SDP. The authors report no duality gap on some standard meshed test systems and randomized versions of these test systems. The zero duality gap property was thus observed experimentally on standard test systems, and further research resulted in sufficient conditions. This is the case, for example, if the objective function is convex and monotonically increasing with the active power generation, and the network has a radial topology \cite{bose2012quadratically, gan2012branch}.  Another approach aiming at global optimality relies on LR \cite{phan2012}, which is further explained in Section~\ref{section:lag-relax}. The author also describes a spatial branch and bound (B\&B) algorithm to close the gap, should it exist one. The ability of both SDP and LR to decrease the optimality gap within a B\&B framework was evaluated in \cite{gopalakrishnan2012global}. If SDP appeared to be computationally more attractive, it showed that it could be very challenging to reach a significant gap reduction within reasonable time limits, even for small test systems.

Multi-period applications related to energy storage are investigated in \cite{gayme2011optimal}, where the SDP relaxation of \cite{lavaei2012} is successfully applied, as their particular application met the conditions for having no duality gap. The authors of \cite{gopalakrishnan2013global} argue that extending \cite{gopalakrishnan2012global} to a multi-period setting yields a SDP too large for current solvers to be solved efficiently and suggest to relax the time-coupling constraints using LR.  However, it ended up being computationally too heavy to make the B\&B approach worthwhile.

Many papers consider the unit commitment problem over an AC network, which is an instance of a multi-period OPF with discrete variables. For instance in \cite{alguacil2000multiperiod}, a generalized Benders decomposition divides the problem in a linear master problem with discrete variables and non-linear multi-period subproblems. Benders cut are generated from the subproblems to tighten the MIP master problem.

\section{Relaxations description}
\label{sec:relax}

We are looking for a computationally affordable relaxation of the problem stated in Section~\ref{sec:problem-statement} that would offer both a narrow optimality gap and a solution close to be feasible.
The main complexity sources of problem (\ref{eq:1})-(\ref{eq:10}) are the discrete decision variables \eqref{eq:2} and the non-convexity of (\ref{eq:8})-(\ref{eq:10}). Furthermore the problem is large scale because of the time-coupling constraints (\ref{eq:5})-(\ref{eq:7}). If the set of constraints (\ref{eq:8})-(\ref{eq:10}) could be addressed independently, finding an optimal solution of (\ref{eq:1})-(\ref{eq:10}) would result in solving less complex subproblems. This decomposition is particularly attractive because:
\begin{itemize}
\item the large time-coupled problem is now a mixed-integer quadratic program (MIQP) or a mixed-integer linear program (MILP) which are much easier to solve than a MINLP of comparable size;
\item every constraint of (\ref{eq:8})-(\ref{eq:10}) only involves period-specific variables and this non-convex program (NLP) can thus be split in $|\mathcal{T}|$ smaller independent problems.
\end{itemize}
However, these two sets of constraints share the power injection variables appearing in (\ref{eq:5})-(\ref{eq:7}) and in the left-hand sides of (\ref{eq:8})-(\ref{eq:9}). Thus some coordination between those subproblems is required to obtain a solution to (\ref{eq:1})-(\ref{eq:10}).

Such a decomposition has already been proposed in \cite{phan2012} for single-period continuous OPFs, where the coordination between the power and voltage subproblems was performed using LR. The extension of this work to the considered problem statement is presented in Section~\ref{section:lag-relax}. In addition, we introduce in Section~\ref{section:nf-relax} a novel flow-based relaxation for this class of multi-period mixed-integer OPFs. The main idea behind this relaxation is that the power flow equations (\ref{eq:8})-(\ref{eq:9}) can be formulated as a network flow with losses.

\subsection{Lagrangian relaxation}
\label{section:lag-relax}

As previously discussed, the author of \cite{phan2012} proposes a Lagrangian Relaxation (LR) scheme
in which the constraints (\ref{eq:8})-(\ref{eq:10}) are dualized.
He proves that this leads to two independent subproblems: a problem involving the active and reactive power injections, and a quadratic problem involving the voltage variables. 
If we apply the same idea to the problem presented in Section \ref{sec:problem-statement}, we obtain the Lagrangian $L$
as
\begin{align*}
L(&\bm{P}, \bm{Q}, \bm{d}, \bm{e}, \bm{f}, \bl, \bg, \ba, \bb)\\
&= f(\bm{P}) + \bm{c}_f \cdot \bm{d}\\
& + \sum_{(t,i) \in \mathcal{N} \times \mathcal{T}} \bl_i^{(t)} \bigg( (\mathbb{M}\bm{P}^{(t)})_i - \sum_{j \in \mathcal{N}(i)} \Big( g_{ij} ({e_i^\tinyt}^2 + {f_i^\tinyt}^2\\
& \hspace{0.5em} - e_i^\tinyt e_j^\tinyt - f_i^\tinyt f_j^\tinyt) + b_{ij} (e_i^\tinyt f_j^\tinyt - f_i^\tinyt e_j^\tinyt) \Big) \bigg)\\
& + \sum_{(t,i) \in \mathcal{N} \times \mathcal{T}} \bg_i^{(t)} \bigg( (\mathbb{M}\bm{Q}^{(t)})_i - \sum_{j \in \mathcal{N}(i)} \Big( b_{ij} (e_i^\tinyt e_j^\tinyt\\ 
& \hspace{0.5em} + f_i^\tinyt f_j^\tinyt - {e_i^\tinyt}^2 - {f_i^\tinyt}^2) + g_{ij} (e_i^\tinyt f_j^\tinyt - f_i^\tinyt e_j^\tinyt) \Big) \bigg)\\
& + \sum_{(t,i) \in \mathcal{N} \times \mathcal{T}}  \ba_i^{(t)} \big( \underline{V}_i^2 - {e_i^\tinyt}^2 - {f_i^\tinyt}^2 \big)\\
&  + \sum_{(t,i) \in \mathcal{N} \times \mathcal{T}}  \bb_i^{(t)} \big({e_i^\tinyt}^2 + {f_i^\tinyt}^2 - \overline{V}_i^2 \big)\\
\end{align*}
where $\bl$, $\bg  \in \mathbb{R}^{|\mathcal{T}||\mathcal{N}|}$ and $\ba$,  $\bb \in \mathbb{R}^{|\mathcal{T}||\mathcal{N}|}_+$ are the Lagrange multipliers for the relaxed constraints.

Any value of the dual function $g$ defined as 
\begin{align}
g(\bl, \bg, \ba, \bb) = \min_{\substack{\bm{P}, \bm{Q} \\ \bm{d}, \bm{e}, \bm{f}}} & & & L(\bm{P}, \bm{Q}, \bm{d}, \bm{e}, \bm{f}, \bl, \bg, \ba, \bb)\label{ldf:1}\\
\text{s.t.}& & & \text{(\ref{eq:2})-(\ref{eq:7})}\label{ldf:2}
\end{align}
provides a lower bound on the optimal value of the original. The Lagrangian dual bound is obtained
by maximizing $g$, which is known to be a concave function.
Still following the approach of \cite{phan2012}, the relaxation is tightened by introducing, $\forall t \in \mathcal{T}$, the constraints
\begin{align}
\sum_{i \in \mathcal{N}} \underline{V}_i^2 \leq \sum_{i \in \mathcal{N}} ({e_i^\tinyt}^2 + {f_i^\tinyt}^2) \leq \sum_{i \in \mathcal{N}} \overline{V}_i^2 \label{sum_vlimits}
\end{align}
If they are redundant in the original problem, they are not in (\ref{ldf:1})-(\ref{ldf:2}) because (\ref{eq:10}) has been relaxed.

More specifically we can rewrite the problem as
\begin{align}
\max_{\substack{\bl, \bg \\ \ba, \bb}} g(\bl, \bg, \ba, \bb)\hspace{-2em}&\label{lr:prob}\\[-1.3em]
= \max_{\substack{\bl, \bg \\ \ba, \bb}} &\Big\{ L_P^{*}(\bl, \bg) +  L_V^{*}(\bl, \bg, \ba, \bb)\nonumber\\[-1.3em]
&+ \hspace{-1em} \sum_{(t,i) \in \mathcal{N} \times \mathcal{T}}  (\ba_i^{(t)}\underline{V}_i^2 - \bb_i^{(t)}\overline{V}_i^2)\Big\}\nonumber
\end{align}

where the power subproblem $L_P(\bl, \bg)$ is defined as
\begin{align*}
L_P^{*}(\bl, \bg) = &\min_{\substack{\bm{P}, \bm{Q} \\ \bm{d}}} &\hspace{-1em}f(\bm{P})&+ \bm{c}_f \cdot \bm{d}\\[-0.9em]
&&&+ \hspace{-1em} \sum_{(t,i) \in \mathcal{N} \times \mathcal{T}} \bl_i^{(t)} (\mathbb{M}\bm{P}^{(t)})_i\\[-0.5em]
&&&+ \hspace{-1em} \sum_{(t,i) \in \mathcal{N} \times \mathcal{T}} \bg_i^{(t)} (\mathbb{M}\bm{Q}^{(t)})_i\\
&\text{s.t.} \hspace{2em} \text{(\ref{eq:2})-(\ref{eq:7})}\hspace{-20em}&&
\end{align*}
and requires solving a MIQP (or MILP). The voltage subproblem $L_V(\bl, \bg, \ba, \bb)$ is on the other hand defined as
\begin{align*}
&L_V^{*}(\bl, \bg, \ba,~\bb)\hspace{-10em}&&\\
&= \sum_{t \in \mathcal{T}} \Big\{ \min_{\bm{e}^{(t)}, \bm{f}^{(t)}} &&- \sum_{i \in \mathcal{N}} \bl_i^{(t)} \sum_{j \in \mathcal{N}(i)} \Big( b_{ij} (e_i^\tinyt f_j^\tinyt - f_i^\tinyt e_j^\tinyt)\\[-0.5em]
&&&+ g_{ij} ({e_i^\tinyt}^2 + {f_i^\tinyt}^2 - e_i^\tinyt e_j^\tinyt - f_i^\tinyt f_j^\tinyt) \Big)\\
&&&- \sum_{i \in \mathcal{N}} \bg_i^{(t)} \sum_{j \in \mathcal{N}(i)} \Big( b_{ij} (e_i^\tinyt e_j^\tinyt + f_i^\tinyt f_j^\tinyt\\[-0.5em]
&&& - {e_i^\tinyt}^2 - {f_i^\tinyt}^2) + g_{ij} (e_i^\tinyt f_j^\tinyt - f_i^\tinyt e_j^\tinyt) \Big)\\
&&& + \sum_{i \in \mathcal{N}} (\bb_i^{(t)}-\ba_i^{(t)}) ({e_i^\tinyt}^2 + {f_i^\tinyt}^2)\\
& \hspace{5em} \text{s.t}. \hspace{2em} (\ref{sum_vlimits}) \hspace{1em} \Big\} \hspace{-10em}&&
\end{align*}
and consists in solving $|\mathcal{T}|$ independent problems that, even though they are non-convex, can be reformulated as trust-region subproblems and solved efficiently in polynomial time.

The convex problem (\ref{lr:prob}) belongs to the class of non-smooth (i.e. non-differentiable) optimization. If subgradient algorithms \cite{boyd2003subgradient} are frequently use to solve these problems, they have shown serious convergence issues for our particular application in the presence of a nonzero duality gap \cite{gopalakrishnan2012global}. For this reason, we suggest to use a bundle method algorithm \cite{feltenmark2000dual} to solve (\ref{lr:prob}).

\subsection{Network flow relaxation}
\label{section:nf-relax}
In the LR scheme presented in Section~\ref{section:lag-relax}, no information on the topology of the network is used in the power subproblem $L_P$. 
Here we present a relaxation that uses the topological information by coupling the original problem with a network flow. 
As the network flow formulation is a linear relaxation of the power flow equations, it does not account for their non-convexities.
In particular it can be observed that in a linear network flow, the total amount of power produced
is equal to the total amount of power consumed, which is rarely the case in our application.
It is therefore important to tighten the formulation by adding some new constraints that 
accounts for these losses in the lines. In particular, we rely on a reformulation-linearization technique (RLT) approach \cite{sherali1998reformulation} 
that yields a convex envelope of the quadratic constraints coming from the power flow.
As a prerequisite for the network flow formulation, we first introduce some notations:
\begin{itemize}
\item $P^{(t)}_{ij}$ is the active power injected in link $(i,j) \in \mathcal{L}$ at bus $i$, positive when power is withdrawn from bus $i$;
\item $Q^{(t)}_{ij}$  is the reactive power injected in link $(i,j) \in \mathcal{L}$ at bus $i$,  positive when power is withdrawn from bus $i$;
\item ${P_{ij}^{loss}}^{(t)}$  is the active power losses in link $(i,j) \in \mathcal{L}$.
\end{itemize}
Using these variables, the conservation of the power flows through links, taking the losses into account, can be written as, $\forall(i,j) \in \mathcal{L}$:
\begin{align}
P^{(t)}_{ij} + P^{(t)}_{ji} &= {P_{ij}^{loss}}^{(t)}\label{eq:conserv1}\\
Q^{(t)}_{ij} + Q^{(t)}_{ji} &= -\frac{b_{ij}}{g_{ij}} {P_{ij}^{loss}}^{(t)} \big( = {Q_{ij}^{loss}}^{(t)} \big)\label{eq:conserv2}
\end{align} 
and the flow conservation at bus $i \in \mathcal{B}$ as:
\begin{align}
(\mathbb{M}~\bm{P}^{(t)})_i = \sum_{j \in \mathcal{N}(i)} P^{(t)}_{ij}\label{eq:conserv3}\\
(\mathbb{M}~\bm{Q}^{(t)})_i = \sum_{j \in \mathcal{N}(i)} Q^{(t)}_{ij}\label{eq:conserv4}
\end{align}
A connection between these flow variables and the voltage variables $\bm{e}$ and $\bm{f}$ is achieved through the following equations:
\begin{gather}
P^{(t)}_{ij} = g_{ij} ({e_i^\tinyt}^2 + {f_i^\tinyt}^2  - e_i^\tinyt e_j^\tinyt - f_i^\tinyt f_j^\tinyt)\nonumber\\+ b_{ij} (e_i^\tinyt f_j^\tinyt - f_i^\tinyt e_j^\tinyt)\label{elec:1}\\
Q^{(t)}_{ij} = b_{ij} (e_i^\tinyt e_j^\tinyt + f_i^\tinyt f_j^\tinyt - {e_i^\tinyt}^2 - {f_i^\tinyt}^2)\nonumber\\+ g_{ij} (e_i^\tinyt f_j^\tinyt - f_i^\tinyt e_j^\tinyt)\label{elec:2}\\
{P_{ij}^{loss}}^{(t)} = g_{ij} ({e_i^\tinyt}^2 + {e_j^\tinyt}^2 + {f_i^\tinyt}^2 + {f_j^\tinyt}^2\nonumber\\- 2 e_i^\tinyt e_j^\tinyt - 2 f_i^\tinyt f_j^\tinyt)\label{elec:3}
\end{gather} 
which are used together with (\ref{eq:conserv1})-(\ref{eq:conserv4}) to obtain a reformulation of the original problem:
\begin{align}
\min_{\substack{\bm{P}, \bm{Q} \\ \bm{d}, \bm{e}, \bm{f}}} & & & f(\bm{P}) + \bm{c}_f \cdot \bm{d}\label{nf}\\
\text{s.t.}& & & (\ref{eq:2})\text{-}(\ref{eq:7})\nonumber\\
& & & \hspace{-0.75em}\forall (t,i) \in \mathcal{T} \times \mathcal{B}:\nonumber\\
& & & (\ref{eq:10}), (\ref{eq:conserv3})\text{-}(\ref{eq:conserv4})\nonumber\\
& & & \hspace{-0.75em}\forall (t,(i,j)) \in \mathcal{T} \times \mathcal{L}:\nonumber\\
& & & (\ref{eq:conserv1})\text{-}(\ref{eq:conserv2}), (\ref{elec:1})\text{-}(\ref{elec:3})\nonumber
\end{align}
This problem is a mixed-integer quadratically constrained quadratic program (MIQCP), which is non-convex just as the original problem. It is important to note that there are redundant constraints in this formulation. For example, removing (\ref{eq:conserv1})-(\ref{eq:conserv2}) and (\ref{elec:3}) would produce an equivalent mathematical program. However, it does not mean that the relaxed counterparts of these constraints will also be redundant. It has indeed been shown in \cite{ruiz2011using} that such redundancy helps generating tighter relaxations.

Such a problem can be relaxed by replacing bilinear (i.e $x_i x_j$) and quadratic (i.e. $x_i^2$) terms by their McCormick envelopes, which can be generated by following the procedure described in Table \ref{rlt_algo}.
\begin{table}[htb]
\small
\begin{leftbar}
\begin{align*}
\text{Let}~x_i \in [l_i, u_i&]~\text{and}~x_j \in [l_j, u_j]\\
\text{then} \hspace{2.3em} x_i x_j &\rightarrow~ w_{ij}\\
\text{with} \hspace{2.8em} w_{ij} &\geq u_i x_j + u_j x_i - u_i u_j\\
w_{ij} &\geq l_i x_j + l_j x_i - l_i l_j\\
w_{ij} &\leq u_i x_j + l_j x_i - u_i l_j\\
w_{ij} &\leq l_i x_j + u_j x_i - l_i u_j
\end{align*}
\caption{Procedure to replace a bilinear term by its convex envelope.}
\label{rlt_algo}
\end{leftbar}
\end{table}
However, before doing so, it is important to observe that such a relaxation converges towards the original problem as the variable domain is getting smaller, i.e. $\max~(x_i x_j - w_{ij})$ converges to zero as ${(\overline{x}_i - \underline{x}_i)}$ and ${(\overline{x}_j - \underline{x}_j)}$ tends to zero too. In other words, the closer the bounds are, the tighter is the relaxation. Unfortunately, the bounds of $\bm{e}$ and $\bm{f}$ are initially quite loose: $e_i^\tinyt$ and $f_i^\tinyt$ belong to $[-\sqrt{\overline{V}_i}, +\sqrt{\overline{V}_i}]$, $\forall (i,t) \in \mathcal{T} \times \mathcal{B}$. In order to tighten the relaxed problem, it would be interesting to refine these bounds given the set $\mathcal{S}$ of feasible solutions of (\ref{eq:1})-(\ref{eq:10}). Because computing such bounds in the original problem would 
result in the same time-complexity as the original problem,
we rely on a subset of period-specific constraints of (\ref{eq:2})-(\ref{eq:10}) to approximate $\mathcal{S}$. For each time period $t \in \mathcal{T}$, some constraints are removed from the original problem to obtain an approximated set $\tilde{\mathcal{S}}_t$ such that $\mathcal{S}_t \subset \tilde{\mathcal{S}}_t$ with $\mathcal{S}_t$ the projection of the original set of feasible solutions to the set of period-$t$-specific variables. In other words, the resulting bounds of $\bm{e}$ and $\bm{f}$ deduced from sets $\tilde{\mathcal{S}}_t$ are guaranteed not to remove any feasible solution from the original problem. In particular, the set $\tilde{\mathcal{S}}_t$ is defined as:
\begin{align*}
\{~(\bm{P}^{(t)},\bm{Q}^{(t)},\bm{e}^{(t)},\bm{f}^{(t)})~|~\text{(\ref{eq:3})-(\ref{eq:4}),(\ref{eq:8})-(\ref{eq:10}) are not violated}~\}
\end{align*}
and finding the upper and lower bounds of a voltage variable $v$ (i.e. $e^\tinyt_i$ or $f^\tinyt_i$, $\forall (i,t) \in \mathcal{B} \times \mathcal{T}$) is equivalent to solving the following problem:
\begin{align}
\overline{v} / \underline{v} = \underset{\substack{\bm{P}^\tinyt, \bm{Q}^\tinyt \\ \bm{e}^\tinyt, \bm{f}^\tinyt}}{\max/\min} & & & v\label{ts:1}\\
\text{s.t.}& & & (\bm{P}^{(t)},\bm{Q}^{(t)},\bm{e}^{(t)},\bm{f}^{(t)}) \in \tilde{\mathcal{S}}_t\label{ts:2}
\end{align}
Even if this problem is much smaller than the original one, it is still non-convex. For this reason, the bounds on $\bm{e}$ and $\bm{f}$ are finally computed by solving an SDP relaxation \cite{anstreicher2012convex} of (\ref{ts:1})-(\ref{ts:2}). These are the bounds used to build the RLT relaxation of (\ref{nf}).

The last tightening step that we perform is to bound
the variables $P^{(t)}_{ij}$, $P^{(t)}_{ji}$, $Q^{(t)}_{ij}$, $Q^{(t)}_{ji}$ and ${P_{ij}^{loss}}^{(t)}$ by solving the SDP relaxation of (\ref{ts:1})-(\ref{ts:2}) with as objective function their expression in equations (\ref{elec:1})-(\ref{elec:3}).

\section{Quantitative analysis}
\label{sec:results}

\subsection{OPF applications}
\label{sec:benchmark-problem}

In order to benchmark the relaxations presented in Section~\ref{sec:relax}, we focus on two applications of the OPF. The first one is the common minimization of generation costs, where we define the cost function $f(\bm{P})$ as
\begin{gather*}
f_{\text{gen}}(\bm{P}) = \sum_{t \in \mathcal{T}} \sum_{g \in \mathcal{G}} \Big( a_g^\tinyt {P_g^\tinyt}^2 + b_g^\tinyt P_g^\tinyt + c_g^\tinyt \Big)
\end{gather*}
with $\mathcal{G}$ the set of generators. In particular, we consider that the generation costs can vary over time. This is modeled by using time-varying parameters $\{a_g^\tinyt, b_g^\tinyt, c_g^\tinyt\}$. In this context, flexible load can be worthwhile to shift the demand when generation costs are low.

The second application is a curtailment minimization and is an extension of the deterministic version of \cite{gemineactive}. In this case, the cost function $f(\bm{P})$ is defined as
\begin{gather*}
f_{\text{curt}}(\bm{P}) = \sum_{t \in \mathcal{T}} \Big[ c_{\text{curt}} \sum_{g \in \mathcal{G}} \big( \overline{P}_g^\tinyt - P_g^\tinyt \big) + c_{\text{losses}} \sum_{d \in \mathcal{D}} P_d^\tinyt \Big]
\end{gather*}
where the first term represents the curtailment costs and the second term expresses the cost of network losses. Such a cost function is representative of the objective of a distribution system operator that operates a network with distributed generators. Flexible loads can be profitable if their consumption is shifted when production from distributed generators is high, e.g. to avoid congestions or over-voltages without relying too much on curtailment. For both applications, the term $\bm{c}_f\cdot\bm{d}$ must be added to the cost function in order to account for availability fees.

\subsection{Implementation details}

The test program is written in C++ and uses several solver libraries. For LR, a continuous relaxation of the original problem is first solved using IPOPT \cite{ipopt} to initialize Lagrange multipliers and solving the non-smooth problem is done with ConicBundle \cite{conic}. The subproblem $L_P$ is solved with MOSEK~\cite{andersen2003implementing} while $L_V$, after being casted into a minimal eigenvalue problem, is addressed using Eigen \cite{eigen}. For the network flow relaxation (NFR), all SDP relaxations as well as the final convex relaxation are solved with MOSEK.

The primal solutions, computed to evaluate the optimality gap of the relaxed solutions, were obtained using SCIP \cite{scip} configured with IPOPT as NLP solver.

\subsection{Instances}
\label{sec:instances}

An instance is defined by a cost function, a network and a number of periods. Table~\ref{tab:networks} presents the different networks used in the test case (if the original test contains shunt admittances, they are ignored).
\begin{table}[htb]
\footnotesize
\begin{center}
\begin{tabular}{|c|c|c|c|c|}
\hline
 & $|\mathcal{B}|$ & $|\mathcal{G}|$ & $|\mathcal{F}|$ & Source\\
 \hline
(A) & 6 & 3 & 3 & \cite{wood2012power}\\
(B) & 9 & 3 & 3 & \cite{zimmerman2011matpower}\\
(C) & 14 & 5 & 4 & \cite{christie2000power}\\
(D) & 6 & 2 & 2 & \cite{gemineactive}\\
\hline
\end{tabular}
\end{center}
\caption{Networks used for the benchmark.\label{tab:networks}}
\end{table}

The cost function $f_{\text{gen}}$ is tested on (A)-(C) and $f_{\text{curt}}$ on (A)-(D). For the curtailment application on networks (A)-(C), one of the generator (the slack bus) is modified to model a connection with another network. The power injection at the corresponding bus can, within some limits, be  either positive or negative.

The test instances are finally generated by considering these 7 (network, cost function) pairs over 4 and 8 periods to obtain a total of 14 instances.

\subsection{Numerical results}
\label{sec:results}

Numerical results on the 14 instances are presented in \mbox{Table~\ref{tab:results}-\subref{tab:results:4}} and \mbox{Table~\ref{tab:results}-\subref{tab:results:8}}. The relative optimality gap is computed as follow:
\begin{gather*}
\text{gap} = \frac{ub^* - lb}{lb}
\end{gather*}
where  $lb$ is the optimal solution of the relaxed problem (i.e. a lower bound) which can vary for every relaxation used  and $ub^*$ is the best primal solution known, and is a fixed number. For each instance, the reported time is the duration of the program before termination, running on a 2.6~GHz processor and limited to a single core. We observe that both relaxations have similar performance for the optimality gap, in the sense that it is almost always within the same order of magnitude. Concerning the running time performance, there is not an approach that outperforms the other as both relaxations show very diverse results.

\begin{table}[tb]
\footnotesize
\centering
\subfloat[Numerical results for $|\mathcal{T}|=4$.]{\label{tab:results:4}
\begin{tabular}{|c||r|r||r|r|}
\hline
{} 		& \multicolumn{2}{c||}{ LR } & \multicolumn{2}{c|}{ NFR } \\
\hline
Case		& gap (\%) 	& time (s)		& gap (\%) & time (s)	 \\
 \hline
(A)$_\text{gen}$ & 2.37 & 203.7 & 4.27 & 11.1 \\
(B)$_\text{gen}$ & 0.00 & 1.2 & 2.24 & 12.7 \\
(C)$_\text{gen}$ & 0.11 & 143.0 & 5.16 & 84.2 \\
(A)$_\text{curt}$ & 79.69 & 45.0 & 225.72 & 16.0 \\
(B)$_\text{curt}$ & 9.07 & 20.1 & 12.53 & 23.5 \\
(C)$_\text{curt}$ & 648.64 & 140.1 & 593.58 & 163.3 \\
(D)$_\text{curt}$ & 60.90 & 40.9 & 60.99 & 11.3 \\
\hline
\end{tabular}
}\\
\subfloat[Numerical results for $|\mathcal{T}|=8$.]{\label{tab:results:8}
\begin{tabular}{|c||r|r||r|r|}
\hline
{} 		& \multicolumn{2}{c||}{ LR } & \multicolumn{2}{c|}{ NFR } \\
\hline
Case		& gap (\%) 	& time (s)		& gap (\%) & time (s)	 \\
 \hline
(A)$_\text{gen}$ & 2.51 & 2905.2 & 4.50 & 38.7 \\
(B)$_\text{gen}$ & 0.00 & 4.1 & 2.20 & 40.7 \\
(C)$_\text{gen}$ & 0.24 & 780.5 & 5.07 & 254.7 \\
(A)$_\text{curt}$ & 124.86 & 83.9 & 255.16 & 82.7 \\
(B)$_\text{curt}$ & 11.90 & 60.9 & 13.22 & 111.0 \\
(C)$_\text{curt}$ & 879.68 & 414.8 & 649.43 & 1207.9 \\
(D)$_\text{curt}$ & 65.10 & 112.5 & 60.09 & 64.1 \\
\hline
\end{tabular}
}\\
\subfloat[Sum of squared infeasibilities of relaxed solutions for constraints (\ref{eq:8})-(\ref{eq:10}).]{\label{tab:infeas}
\begin{tabular}{|c||r|r||r|r|}
\hline
{} 		& \multicolumn{2}{c||}{ $|\mathcal{T}|=4$ } & \multicolumn{2}{c|}{ $|\mathcal{T}|=8$ } \\
\hline
Case		& \multicolumn{1}{c|}{LR} 	& \multicolumn{1}{c||}{NFR} & \multicolumn{1}{c|}{LR} 	& \multicolumn{1}{c|}{NFR}  \\
 \hline
(A)$_\text{gen}$ & 6.02 & 0.02 & 8.72 & 0.05 \\
(B)$_\text{gen}$ & 70.87 & 99.33 & 141.50 & 196.21\\
(C)$_\text{gen}$ & 1.24 & 1.72 & 1.86 & 3.58\\
(A)$_\text{curt}$ & 86.75 & 6.78 & 163.28 & 13.46\\
(B)$_\text{curt}$ & 179.86 & 152.51 & 142.72 & 162.40\\
(C)$_\text{curt}$ & 456.88 & 6.91 & 57.57 & 16.35\\
(D)$_\text{curt}$ & 854.16 & 0.10 & 1564.37 & 0.19\\
\hline
\end{tabular}
}\\
\caption{Results for the 14 instances.\label{tab:results}}
\end{table}

We are also interested in evaluating another feature of these relaxations: the level of infeasibility of their solutions in the original problem. This feature can indeed affect the efficiency of a relaxation within a spatial B\&B framework \cite{lawler1966branch} when seeking for a globally optimal solution of \mbox{Problem~(\ref{eq:1})-(\ref{eq:10})}. Relaxed solutions that are closer to feasibility can speed up the discovery of feasible solutions and at the same time provide upper bounds to the objective function earlier in the space exploration procedure. Obtaining upper bounds is critical for these approaches as it helps pruning nodes and reduces the computational budget required before termination. \mbox{Table~\ref{tab:results}-\subref{tab:infeas}} presents the sum of squared infeasibilities for the set of constraints (\ref{eq:8})-(\ref{eq:10}) (i.e. those relaxed in LR and NFR). We observe that NFR shows less infeasibility than LR on 9 out of 14 instances. For some cases, NFR produces solutions that are very close to be feasible (e.g. (A)$_\text{gen}$ and (D)$_\text{curt}$) while LR does not exhibit similar performances even when it is able to close the gap (e.g. (B)$_\text{gen}$). In addition, some of the solutions of LR are affected with a very high level of infeasibility (e.g. (C)$_\text{curt}$ and (D)$_\text{curt}$), which is orders of magnitude worse than NFR.


\section{Conclusion}
\label{sec:conclusion}

In this paper, we present a novel relaxation for multi-period OPF with discrete variables that is based on a network-flow reformulation. While the lower bounds it produces are comparable with the Lagrangian relaxation, the infeasibility of the relaxed solutions is reduced. This feature suggests that it is worthwhile to evaluate NFR beside the current state-of-the-art relaxations (i.e.  \cite{phan2012} and \cite{lavaei2012}) within  a B\&B framework.

On the other hand, this relaxation should still be improved on two aspects. The first one is the quality of lower bounds, especially for curtailment applications.  We believe that a special care should be taken concerning the upper bounds of the active losses in links. We observed that the SDP relaxation used to compute these bounds is not very informative and it penalizes the tightness of the overall relaxation. The second aspect to improve is on the computational side. For this purpose, we would like to consider subnetworks instead of the whole network to infer the bounds on the voltage and link-flow variables. If it would reduce the size of SDP problems and speed up their convergence, it could also reduce the value of the resulting bounds. For this reason, an iterative approach that would increase the size of specific subproblems to narrow the most useful bounds is not to put aside. 

Following the observations of this work, we think that another interesting research direction would be to merge the two relaxations considered in this paper. Tightening the power subproblem of a Lagrangian relaxation with a network-flow relaxation could both improve the convergence of the non-smooth problem of LR thanks to a tighter subproblem and reduce the infeasibility of produced solutions.

\section*{Acknowledgments}
This research is supported by the public service of Wallonia – Department of Energy and Sustainable Building within the framework of the GREDOR project.
The authors thank the financial support of the Belgian Network DYSCO, an Inter-university Attraction Poles Program initiated by the Belgian State, Science Policy Office.

\printbibliography

\end{document}